\def\IK{{\Bbb K}}

\def\IC{\Bbb C} 
\def\ID{{\Bbb D}}

\def\zbar{{\overline{z}}}

\documentclass[12pt]{article} 

\usepackage[psamsfonts]{amssymb} 
\usepackage{amsfonts,amsmath}

\usepackage{epsfig,multicol}

\newtheorem{theorem}{Theorem}%[section] 
\newtheorem{lemma}{Lemma}%[section] 
\title{Non-variational extrema of exponential Teichm\"uller spaces}
\date{}
\author{Gaven Martin \and Cong Yao \thanks{Both author's research supported in part by the New Zealand Marsden Fund\newline \newline This work forms part of the second authors PhD Thesis.} }  \date{{\small to Pekka Koskela on the ocassion of his $60^{th}$ birthday} }
\begin{document}
\maketitle  
\begin{abstract}  The exponential Teichm\"uller spaces $E_p$,  $0\leq p \leq \infty$, interpolate between the classical Teichm\"uller space ($p=\infty$) and the space of harmonic diffeomorphisms $(p=0)$.  In this article we prove the existence of non-variational critical points for the associated functional:  mappings $f$ of the disk whose distortion is $p$-exponentially integrable, $0<p<\infty$,  yet for {\em any} diffeomorphism $g(z)$ of $\ID$   with $g|\partial\ID=identity$ and $g\neq identity$ we have $f\circ g$ is not of $p$-exponentially integrable distortion.
\end{abstract}
 \section{Introduction}  In 1939 Teichm\"uller stated his famous theorem \cite{T}: In the homotopy class of a diffeomorphism between closed Riemann surfaces, there is a unique extremal quasiconformal mapping of smallest maximal distortion $\IK(z,f)$. Furthermore, this function $f$ is either conformal or has a Beltrami coefficient $\mu_f=f_\zbar/f_z$ of the form
\begin{equation}\label{teich1}
\mu_f(z) =k\; \frac{\overline{\Psi(z)}}{|\Psi(z)|},
\end{equation}
where $\Psi$ is a holomorphic function and $0\leq k <1$ is a constant.  Indeed $\Psi$ defines a holomorphic quadratic differential providing local coordinates in which the extremal map is linear.  In \cite{Ah3} Ahlfors gave the first complete proof of Teichm\"uller's theorem using a variational approach essentially minimising the $L^p$-norm of the distortion and letting $p\to\infty$.  What is remarkable here is that in Ahlfors' approach the $p$-minimisers have essentially no regularity or good topological properties - they may not even be continuous (and 70 years later we are just beginning to establish some of these facts) but the associated holomorphic Hopf differentials converge as $p\to\infty$ to yield the Beltrami coefficient in the limit.

In our work \cite{MY},  following initial explorations in \cite{IMO},  we developed the exponential Teichm\"uller spaces where instead of minimising the maximal distortion $\IK(z,f)$ or its $L^p$-norm we minimize the $p$-exponential norm 
\begin{equation}\label{2}
{\cal E}_p[f] =\int_\Sigma e^{p\IK(z,f)}\; dz
\end{equation}  
where $f$ is in a given homotopy class or $f$ has fixed boundary data. 

Due to the many recent and significant advancements in the theory of mappings of finite distortion (see \cite{AIM} as a starting point) we were able to show that for each $p$ a homeomorphic minimiser $f_p$ exists,  that as $p\to\infty$ we obtain the extremal quasiconformal mapping (giving a new proof for Teichm\"uller's theorem) as a limit $f_p\to f_\infty$ and that as $p\to 0$ we get the harmonic mappings as a limit.  This ``unified'' the classical Teichm\"uller theory and the harmonic mapping approach to Teichm\"uller theory of \cite{Wentworth,Tromba,Wolf}.   Further, if the distortion of this minimiser is locally bounded then the minimiser is a diffeomorphism.  Actually only a much weaker condition is necessary: a the mapping should satisfy the (inner) variational equations and have distortion in a slightly better exponential class.  We conjecture all minimisers have this slight additional regularity,  but we know that some critical points do not. It seems this is necessary to get the strong regularity of critical points.   

\medskip

In this article we construct non-variational critical points $f\mapsto {\cal E}_p[f]$ for the functional defined at (\ref{2}) on the disk, so $\Sigma=\ID$ and $f|\partial\ID=f_0$ a given homeomorphism. Lifting to the universal cover will give the result on arbitrary Riemann surfaces.  In particular our main result here is the following.

\begin{theorem}\label{mainthm} Let $0<p<\infty$.  Then there is a homeomorphism $f:\overline{\ID}\to\overline{\ID}$ such that 
\[{\cal E}_p[f]=\int_\ID e^{p\IK(z,f)}\; dz < \infty, \]
$f\in W^{1,2}_{loc}(\ID)$,  the Sobolev space of mappings with locally square integrable first derivatives, and which has the following property:  If $g:\ID\to\ID$  is a diffeomorphism which is not conformal,  then 
\[ {\cal E}_p[f\circ g] = +\infty  \]
\end{theorem}
For the variational problem we would want the diffeomorphism to be the identity on the boundary,  but that is not necessary here.  In fact our solution is entirely local.  If $U\subset\ID$ is open and $g|U$ is not conformal,  then we will prove that 
\[\int_U e^{p\IK(z,f\circ g)}\; dz = +\infty \]
The reader will easily see that this is covered in the reduction to the following lemma and its proof,  however we do make some comments.  We set the lemma up in a variational fashion to be more useful in other applications.  The restriction $0<|t|<1$ and $\|\nabla \varphi\|<1$ are only to ensure the variation is through diffeomorphisms.

\begin{lemma}\label{mainlem} Let $0<p<\infty$.  Then there is a homeomorphism $f:\overline{\ID}\to\overline{\ID}$ such that 
\[{\cal E}_p[f]=\int_\ID e^{p\IK(z,f)}\; dz < \infty, \]
$f\in W^{1,2}_{loc}(\ID)$,  the Sobolev space of mappings with locally square integrable first derivatives, and which has the following property.  Let $\varphi\in C^{\infty}_{0}(\ID)$ be a non-zero smooth test function with $\|\nabla \varphi\|_{L^\infty(\ID)}<1$. For all $0<t<1$ set $f^t(z) = f(z+t\varphi(z))$.  Then $f^t:\ID\to\ID$ is a homeomorphism with $f^t|\partial\ID=f|\partial\ID$,  $f^t\in W^{1,2}_{loc}(\ID)$,  $E_{q}[f^t]<\infty$ for all $q\leq q_t <q$ with $q_t\to p$ as $t\to 0$,  However
\[ {\cal E}_p[f^t] = +\infty, \hskip20pt 0<t<1. \]
\end{lemma}

\section{Proof of Lemma \ref{mainlem}.}

\subsubsection{A basic example.}
\noindent We start with the function
\begin{equation}\label{Fr}
F(r)=\frac{1}{r^2\log^2(\frac{2}{r})},\quad0<r<1.
\end{equation}
This function satisfies
\[
\int_0^1F(r)rdr<\infty,
\hskip5pt \mbox{and for any $q>1$},\hskip5pt
\int_0^1F^q(r)rdr=+\infty.
\]
The indefinite integral can be evaluated in terms of special functions,  but with a bit of work one can find
\begin{equation}\label{Fint} \int_0^1F(r) \, r\; dr  = \frac{1}{\log 2}.\end{equation}
We would therefore like
\begin{equation}
e^{\IK(z,f)}=\frac{e}{r^2\log^2(\frac{2}{r})},\quad r=|z|.
\end{equation}
and so $
\IK(z,f)=1-2\log(r\log\frac{2}{r}). $
Given an increasing surjective function $\rho:[0,1]\to[0,1]$,  we can define the {\em radial stretching}
\[
f(z)=\frac{z}{|z|}\rho(|z|).
\]
From \cite[Section 2.6]{AIM}) we  compute the Beltrami coefficient
\begin{equation}
\mu_f(z)=\frac{z}\zbar \frac{|z|\dot\rho(|z|)-\rho(|z|)}{|z|\dot\rho(|z|)+\rho(|z|)}.
\end{equation}
Then,
\[
\IK(z,f)=\frac{1+|\mu(z)|^2}{1-|\mu(z)|^2}=\frac{(r\dot\rho+\rho)^2+(r\dot\rho-\rho)^2}{(r\dot\rho+\rho)^2-(r\dot\rho-\rho)^2}=\frac{1}{2}(\frac{r\dot\rho}{\rho}+\frac{\rho}{r\dot\rho}).
\]
Equivalently
\[
\frac{r\dot\rho}{\rho}=\IK(z,f)\pm\sqrt{\IK^2(z,f)-1}.
\]
Choosing the larger answer, which is no smaller than 1,  we have
\begin{eqnarray*}
\frac{r\dot\rho}{\rho} & = & 1-2\log(r\log\frac{2}{r})+\sqrt{[1-2\log(r\log\frac{2}{r})]^2-1}, \nonumber \\
\log\rho(r) & = & \int_1^r\frac{1-2\log(s\log\frac{2}{s})+\sqrt{[1-2\log(s\log\frac{2}{s})]^2-1}}{s}ds,\nonumber \\
\rho(r)& = & \exp[\int_1^r\frac{1-2\log(s\log\frac{2}{s})+\sqrt{[1-2\log(s\log\frac{2}{s})]^2-1}}{s}ds].
\end{eqnarray*}
Note this satisfies $\rho(1)=e^0=1$, and $\log\rho(0)\leq\int_1^0\frac{1}{s}ds=-\infty$, so $\rho(0)=0$ as required. This gives us a self-homeomorphism of $\ID $
\begin{equation}\label{fdef}
f(z)=\frac{z}{|z|}\rho(|z|),
\end{equation}
whose distortion $\IK(z,f)$ is exponentially integrable,  ${\cal E}_1[f]<\infty$.  As $f$ is a homeomorphism with exponentially integrable distortion,  it extends homeomorphically to the boundary circle.

Next we recall the composition formula.  Given $g^t:\ID\to\ID$,
\[
\IK(g^t(z),f\circ(g^t)^{-1})=\IK(z,f)\IK(z,g^t)[1-\frac{4\Re e(\mu_f\overline{\mu_{g^t}})}{(1+|\mu_f|^2)(1+|\mu_{g^t}|^2)}].
\]
We are interested in the sign of the term $\Re e(\mu_f\overline\mu_{g^t})$. If it is non-positive, we then have
\[
\IK(g^t(z),f\circ(g^t)^{-1})\geq\IK(z,f)\IK(z,g^t).
\]
Observe that for each pair of complex numbers $z=a+bi$ and $w=c+di$, we have
\[
\Re e(z\bar{w})=ac+bd.
\]
Thus, at least one of the following is non-positive:
\[
\Re e(z\bar{w}),\quad\Re e(zw),\quad\Re e(-z\bar{w}),\quad\Re e(-zw).
\]
Next,  for a radial stretching  $f(z)=\frac{z}{|z|}\rho(|z|)$,
\begin{eqnarray*}
\mu_f(z)=\frac{z}\zbar \frac{|z|\dot\rho(|z|)-\rho(|z|)}{|z|\dot\rho(|z|)+\rho(|z|)}, &&
\mu_f(iz)=\frac{iz}{\overline{iz}}\frac{|z|\dot\rho(|z|)-\rho(|z|)}{|z|\dot\rho(|z|)+\rho(|z|)}=-\mu_f(z),\\
\mu_f(\bar{z})=\frac\zbar {z}\frac{|z|\dot\rho(|z|)-\rho(|z|)}{|z|\dot\rho(|z|)+\rho(|z|)}=\overline{\mu_f(z)}, &&
\mu_f(i\bar{z})=\frac{i\bar{z}}{\overline{i\bar{z}}}\frac{|z|\dot\rho(|z|)-\rho(|z|)}{|z|\dot\rho(|z|)+\rho(|z|)}=-\overline{\mu_f(z)}.
\end{eqnarray*}
Given $\varphi\in C^{\infty}_{0}(\ID)$ with $\|\nabla\varphi\|_{L^\infty(\ID)} <1$ the mapping $g^t(z)=z+t \varphi(z)$ is a diffeomorphism of $\ID$ which is the identity in a neighbourhood of $\partial\ID$. For $|t|<1$ we compute that 
\[ \mu_{g^t}=\frac{t\varphi_\zbar }{1+t\varphi_z} \]
 is smooth and compactly supported on $\ID$. 
 
 Assume that $\Re e(\mu_{g^t}(0))\neq0$ and $\Im m(\mu_{g^t}(0))\neq0$. By choosing $r^\prime$ small enough, there is a neighbourhood $A=\ID (0,r^\prime)$ in which $\Re e(\mu_{g^t})$ and $\Im m(\mu_{g^t})$ do not change their signs,  and further $
\inf_{z\in A}|\mu_{g^t}(z)|\geq\varepsilon_1>0. $
Then,
\[
q:=\inf_{z\in A}\IK(z,g^t)>1.
\]
Also, $J(z,g^t)=|1+t\phi_z|^2-t^2|\phi_\zbar |^2$ shows us  that $
J(z,g^t)\geq c_t>0$.

Combining these facts we conclude
\begin{eqnarray*}
\lefteqn{\int_{\ID }\exp[\IK(w,f\circ(g^t)^{-1})]dw}\\ &=& \int_\ID \exp[\IK(z,f)\IK(z,g^t)-\frac{4\Re e(\mu_f\overline\mu_{g^t})}{(1-|\mu_f|^2)(1-|\mu_{g^t}|^2)}]J(z,g^t)dz\\
&\geq&c_t\int_A\exp[\IK(z,f)\cdot\inf_{z\in A}\IK(z,g^t)]dz\\
&=&c_t \int_A[\frac{e}{|z|^2\log^2(\frac{2}{|z|})}]^qdz = c_t\int_0^{2\pi}d\theta\int_0^{r^\prime}[\frac{e}{r^2\log^2(\frac{2}{r})}]^q\cdot rdr\\
&\geq& 2\pi c_t \int_0^{r^\prime}[\frac{e}{r\log^2(\frac{2}{r})}]^qdr = +\infty.
\end{eqnarray*}
We now consider the condition $\Re e(\mu_{g^t}(0))\Im m(\mu_{g^t}(0))\neq0$. Since $\mu_{g^t}(z)=\frac{t\varphi_\zbar (z)}{1+t\varphi_z(z)}$ we can compute
\[
\Re e(\mu_{g^t})=\frac{1}{2}\Big(\frac{t\varphi_\zbar }{1+t\varphi_z}+\frac{t\overline{\varphi_\zbar }}{\overline{1+t\varphi_z}}\Big)=\frac{t\Re e(\varphi_\zbar )+t^2\Re e(\varphi_\zbar \overline{\varphi_z})}{|1+t\varphi_z|^2},
\]
\[
\Im m(\mu_{g^t})=\frac{1}{2i}\Big(\frac{t\varphi_\zbar }{1+t\varphi_z}-\frac{t\overline{\varphi_\zbar }}{\overline{1+t\varphi_z}}\Big)=\frac{t\Im m(\varphi_\zbar )+t^2\Im m(\varphi_\zbar \overline{\varphi_z})}{|1+t\varphi_z|^2}.
\]
So the condition $\Re e(\mu_{g^t}(0))\Im m(\mu_{g^t}(0))\neq0$ is satisfied only if
\begin{equation}
\Re e(\varphi_\zbar (0))\Im m(\varphi_\zbar (0))\neq 0.
\end{equation}
We record this as follows.
\begin{theorem}
The Sobolev homeomorphism $f$ defined at (\ref{fdef}) has exponentially integrable distortion $\IK(z,f)$. However, for any $\varphi\in C_0^\infty(\ID )$ with 
\begin{itemize}
\item $\|\nabla\varphi\|<1$,
\item $\Re e(\varphi_\zbar (0))\Im m(\varphi_\zbar (0))\neq 0$,
\item $g^t$ defined as above,
\end{itemize} then for any $t\in(-1,1)$, $t\neq 0$,
\[
\int_{\ID }\exp[\IK(w,f\circ(g^t)^{-1})]dw=+\infty.
\]
\end{theorem}

\subsection{Extending the basic example.}
\noindent To complete the proof of the lemma,  and thus of the main theorem, we have to remove the condition $\Re e(\varphi_\zbar (0))\Im m(\varphi_\zbar (0))\neq 0$.  We will have to modify $f$ to do so.  To this end we need the following lemma:
\begin{lemma}\label{lemma1}
In the unit disk $\ID $, there is a countable dense subset $\{z_k\}$, disjoint Borel sets $S_1$, $S_2$, $S_3$, $S_4$, and a positive number $\delta>0$, such that at every point $z_k$, there is an $R_k>0$ that for any $0<r<R_k$ and $i=1,2,3,4$,
\[
\frac{|S_i\cap D(z_k,r)|}{|D(z_k,r)|}>\delta.
\]
\end{lemma}
We postpone the proof of this lemma to the next section. Now let $F(r)$ be as defined in (\ref{Fr}). This time we define the distortion function as
\[
\IK(z,f)=1+\frac{1}{p}\log\Big[\sum_k\frac{1}{2^k}F(|z-z_k|)\chi_{D(z_k,dist(z_k,\partial\ID ))}\Big],
\]
where $\{z_k\}\subset\ID $ is a dense subset as in Lemma \ref{lemma1}. Then $\exp[p\IK(z,f)]\in L^1(\ID )$. Indeed, using (\ref{Fint}) we see
\[
\int_\ID \exp[p\IK(z,f)]dz\leq\sum_k\frac{2\pi e^p}{2^k}\int_0^1F(r)rdr = \frac{2\pi e^p}{\log 2} <\infty.
\]
The absolute value of the Beltrami coefficient that corresponds to $\IK(z,f)$ is
\begin{equation}
|\mu_f(z)|=\sqrt{\frac{\IK(z,f)-1}{\IK(z,f)+1}}.
\end{equation}
By the existence theorem for mappings with exponentially integrable distortion,  \cite[Theorem 20.4.9]{AIM},  we may set $\mu_f(z)=|\mu_f(z)|e^{i\theta(z)}$ for any measurable function $\theta:\ID \to[0,2\pi)$ and then find a homeomorphism $f:\overline{\ID }\to\overline{\ID }$ that has Beltrami coefficient $\mu_f$ on $\ID $. As above,  using the  composition formula 
\begin{equation}
\IK(g^t(z),f\circ(g^t)^{-1})=\IK(z,f)\IK(z,g^t)[1-\frac{4\Re e(\mu_f\overline\mu_{g^t})}{(1+|\mu_f|^2)(1+|\mu_{g^t}|^2)}].
\end{equation}
\begin{equation}
\Re e(\mu_f\overline{\mu}_{g^t})=\frac{1}{|1+t\varphi_z|^2}\Big(t\Re e(\varphi_\zbar \overline{\mu_f})+t^2\Re e(\varphi_\zbar \overline{\mu_f\varphi_z})\Big).
\end{equation}
We now determine the argument of $\mu_f$. Let $S_1$, $S_2$, $S_3$, $S_4$ as in Lemma 2.3.2, and set
\[
\mu_f(z)=\begin{cases}
|\mu_f(z)|,& z\in S_1;\\
-|\mu_f(z)|,& z\in S_2;\\
i|\mu_f(z)|,& z\in S_3;\\
-i|\mu_f(z)|,& z\in\ID -\bigcup_{i=1}^3S_i\supset S_4.
\end{cases}
\]
So each has density $\delta$ near every $z_i$. We put $\pm|\mu_f|$ and $\pm i|\mu_f|$ respectively  into the above and get
\begin{equation}
\Re e(\mu_f\overline{\mu}_{g^t})=\frac{|\mu_f|}{|1+t\varphi_z|^2}\Big(t\Re e(\varphi_\zbar )+t^2\Re e(\varphi_\zbar \overline{\varphi_z})\Big),\quad z\in S_1;
\end{equation}
\begin{equation}
\Re e(\mu_f\overline{\mu}_{g^t})=-\frac{|\mu_f|}{|1+t\varphi_z|^2}\Big(t\Re e(\varphi_\zbar )+t^2\Re e(\varphi_\zbar \overline{\varphi_z})\Big),\quad z\in S_2;
\end{equation}
\begin{equation}
\Re e(\mu_f\overline{\mu}_{g^t})=\frac{|\mu_f|}{|1+t\varphi_z|^2}\Big(t\Im m(\varphi_\zbar )+t^2\Im m(\varphi_\zbar \overline{\varphi_z})\Big),\quad z\in S_3;
\end{equation}
\begin{equation}
\Re e(\mu_f\overline{\mu}_{g^t})=-\frac{|\mu_f|}{|1+t\varphi_z|^2}\Big(t\Im m(\varphi_\zbar )+t^2\Im m(\varphi_\zbar \overline{\varphi_z})\Big),\quad z\in S_4,
\end{equation}
and we recall here that $0<|t|<1$.  There are now three cases to consider depending on $\varphi$:\\

\noindent{\bf (1)} $\varphi_\zbar \equiv0$ in $\ID $. Since $\varphi\in C_0^\infty(\ID )$, this happens only when $\varphi\equiv0$ in $\ID $, and then $g^t(z)=z$, $f\circ(g^t)^{-1}=f$.\\

\noindent{\bf (2)}  Suppose $\Re e(\varphi_\zbar )$ is not the constant $0$, say $\Re e(\varphi_\zbar )(z_0)>0$ at some point $z_0\in\ID $. Then by the smoothness of $\varphi$ there is an open neighbourhood $U$ where $\Re e(\varphi_\zbar )\geq\varepsilon_1>0$. Following the basic example,   for any $t\in(-1,0)$ in $U\cap S_1$ we have
$|\mu_{g^t}|\geq\varepsilon_2t>0$, $
\Re e(\mu_f\overline{\mu}_{g^t})<0,$
and
$J(z,g^t)=|1+t\phi_z|^2-t^2|\phi_\zbar |^2>c_t.$
Then,  as before, 
\[
\IK(g^t(z),f\circ(g^t)^{-1})\geq q\IK(z,f).
\]
The density of $\{z_k\}$ implies there is  $z_k\in U$. We choose a small disk $D(z_k,r_0)\subset A$, where $r_0<R_k$ as in Lemma \ref{lemma1}. Following our earlier arguments   we now can compute
\begin{eqnarray*}
\lefteqn{\int_\ID \exp[p\IK(w,f\circ(g^t)^{-1})]dw}\\&=&\int_\ID \exp[p\IK(g^t(z),f\circ(g^t)^{-1})J(z,g^t)]dz \geq\frac{1}{2}\int_{D(z_k,r_0)\cap S_1}\exp[pq\IK(z,f)]dz\\
&\geq& \frac{\delta}{2}\int_{D(z_k,r_0)}\exp[pq\IK(z,f)]dz \geq\pi\delta\int_0^{r_0}(\frac{e^p}{2^k}F(r))^qrdr\\
&\geq & C\int_0^{r_0}F^q(r)rdr=\infty.
\end{eqnarray*}

If $\Re e(\varphi_\zbar )(z_0)<0$, then of course the same result follows.\\

\noindent{\bf (3)} Suppose $\Im m(\varphi_\zbar )$ is not the constant $0$. Then the result follows in an entirely similar fashion.\\

This last observation completes the proof. \hfill $\Box$

\medskip

 In fact we can  generalise Lemma \ref{mainlem} to the complex coefficient case. That is, set
\[
g^\eta=z+\eta\varphi,\quad\eta\in\IC ,\quad\varphi\in C_0^\infty(\ID ).
\]
In this case
\begin{equation}
\Re e(\mu_f\overline{\mu}_{g^\eta})=\frac{1}{|1+\eta\varphi_z|^2}\Big(\Re e(\eta\varphi_\zbar \overline{\mu_f})+\Re e(\eta^2\varphi_\zbar \overline{\mu_f\varphi_z})\Big).
\end{equation}
Again we choose $\mu_f$ and $S_1$, $S_2$, $S_3$, $S_4$ same as before. Write $\eta=|\eta|e^{i\alpha}$, then  
\[
\Re e(\mu_f\overline{\mu}_{g_\eta})=\frac{|\mu_f|}{|1+\eta\varphi_z|^2}\Big(|\eta|\Re e(e^{i\alpha}\varphi_\zbar )+|\eta|^2\Re e(e^{2i\alpha}\varphi_\zbar \overline{\varphi_z})\Big),\quad z\in S_1,
\]
and analogously for $S_2$, $S_3$, $S_4$.\\

Then, for any non-constant $\varphi$, we may find a neighbourhood in $\ID $ where either $\Re e(e^{i\alpha}\varphi_\zbar )$ or $\Im m(e^{i\alpha}\varphi_\zbar )$ is nonzero. So by the same argument as above
\[
\int_\ID \exp[p\IK(g_\eta(z),f\circ(g_\eta)^{-1})]dz=\infty.
\]
For every fixed $\varphi$, the number $\varepsilon$ depends only on $\alpha$. If we let $\alpha$ move on $[0,2\pi]$, then $e^{i\alpha}\varphi_\zbar $ and $e^{2i\alpha}\varphi_\zbar \overline{\varphi_z}$ move continuously  and thus $\varepsilon=\varepsilon(\alpha)$ can be chosen as a continuous function of $\alpha$. Since $[0,2\pi]$ is compact, $\varepsilon(\alpha)$ admits a positive minimal value. Then we have established the following theorem.

\begin{theorem}\label{local}
For every $p>0$, there is a homeomorphism $f:\ID \to\ID $ such that $\exp[p\IK(z,f)]\in L^1(\ID )$, and for any non-constant $\varphi\in C_0^\infty(\ID )$  and for any $\eta\in\IC\setminus\{0\} $, $\exp[p\IK(g_\eta(z),f\circ(g_\eta)^{-1})]\notin L^1(\ID )$, where $g_{\eta}=z+\eta\varphi$.
\end{theorem}
\noindent{\bf Proof.} There is really only one thing to note,  and that is that the  gradient estimate we used above,  namely  $\|\nabla \varphi\|_{L^{\infty}(\ID)}<1$, and the condition $t\in [-1,1]$ were only used to ensure that $g^\eta$ is a homeomorphism.  The do not affect the divergence of the integrand  $\exp[p\IK(g_\eta(z),f\circ(g_\eta)^{-1})]$ in $L^1(\ID)$. \hfill $\Box$

\medskip

As an easy consequence,  the result claimed in Theorem \ref{mainthm} now follows directly.  This is because Theorem \ref{local} is entirely local and any diffeomorphism $g:\ID\to\ID$ which is the identify on the boundary of the disk,  and not equal to the identity, must have a point $z_0\in \ID$ where the distortion is greater than 1.  Choose a $C^{\infty}_{0}(\ID)$ function $\varphi$,  $\varphi(z)=1$ near $z_0$ and note that  
\[  g(z) = z+ \varphi(z)(g(z)-z)   \]
near $z_0$.  The result then follows from Theorem \ref{local}. \hfill $\Box$

\subsubsection{The proof of Lemma \ref{lemma1}}
\noindent As a first step, we start with the point $p^1=(0,0)$ and choose the disk sectors
\[
S_1^1=\{z\in D(p^1,\frac{1}{2^4}):0<\arg(z-p^1)<\frac{\pi}{2}\},
\]
\[
S_2^1=\{z\in D(p^1,\frac{1}{2^4}):\frac{\pi}{2}<\arg(z-p^1)<\pi\},
\]
\[
S_3^1=\{z\in D(p^1,\frac{1}{2^4}):\pi<\arg(z-p^1)<\frac{3\pi}{2}\},
\]
\[
S_4^1=\{z\in D(p^1,\frac{1}{2^4}):\frac{3\pi}{2}<\arg(z-p^1)<2\pi\}.
\]
We construct inductively. At Step $n\geq2$, we choose the points $p^n_{j,l}=(\frac{j}{2^{n-1}},\frac{l}{2^{n-1}})$, for any integers $j,l\in[1-2^{n-1},2^{n-1}-1]$ such that $D(p^n_{j,l},\frac{1}{2^{4n}})\subset\ID $, and $p^n_{j,l}$ has not been chosen in the previous steps. Define the sector unions
\[
S_1^n=\bigcup_{j,l}\{z\in D(p^n_{j,l},\frac{1}{2^{4n}}):0<\arg(z-p^n_{j,l})<\frac{\pi}{2}\},
\]
\[
S_2^n=\bigcup_{j,l}\{z\in D(p^n_{j,l},\frac{1}{2^{4n}}):\frac{\pi}{2}<\arg(z-p^n_{j,l})<\pi\},
\]
\[
S_3^n=\bigcup_{j,l}\{z\in D(p^n_{j,l},\frac{1}{2^{4n}}):\pi<\arg(z-p^n_{j,l})<\frac{3\pi}{2}\},
\]
\[
S_4^n=\bigcup_{j,l}\{z\in D(p^n_{j,l},\frac{1}{2^{4n}}):\frac{3\pi}{2}<\arg(z-p^n_{j,l})<2\pi\}.
\]
Write
\[
S^n=S_1^n\cup S_2^n\cup S_3^n\cup S_4^n.
\]
We now define $S_i$ as the set that $z\in S_i^n$ for some $n$ but not in $S^m$ for any $m\geq n+1$. Precisely,
\[
S_i=\bigcup_{n=1}^\infty\Big(S_i^n\cap\bigcap_{m=n+1}^\infty(S^m)^c\Big),\quad i=1,2,3,4.
\]
We claim that the points $p^1$, $p^n_{j,l}$ and the sets $S_i$ satisfy the requirements.\\

We estimate the total area of $\bigcup_{n\geq1} S^n$. At each Step $n$, we have no more than $2^{2n}$ points, and each disk has area $\frac{\pi}{2^{8n}}$. Thus
\begin{equation}
\Big|\bigcup_{n=1}^\infty S^n\Big|\leq\sum_{n=1}^\infty\frac{2^{2n}}{2^{8n}}\pi<\frac{\pi}{2^5}.
\end{equation}
Fix any point $p=p^n_{j,l}$. Then
\[
\{z\in D(p,\frac{1}{2^{4n}}):0<\arg(z-p)<\frac{\pi}{2}\}\subset S_1^n.
\]
Let $r<\frac{1}{2^{4n}}$ be an arbitrary number. Let $N$ be the largest integer that $\frac{1}{2^{N}}>r$. Consider the sector
\[
F:=\{z\in D(p,r):0<\arg(z-p)<\frac{\pi}{2}\}\subset S_1^n.
\]
Note that by the choice of $N$, $F\cap S^m=\emptyset$ for any $n+1\leq m\leq N-1$. So we consider the disk $D(p,\frac{1}{2^N})$. Analogously  we have
\[
|D(p,\frac{1}{2^N})\cap\bigcup_{n\geq N}S^n|\leq\frac{1}{2^{2N}}\cdot\frac{\pi}{2^5}.
\]
On the other hand, by the choice of $N$ we have $\frac{1}{2^{N+1}}\leq r<\frac{1}{2^N}$. So
\[
|F|\geq\frac{1}{4}\cdot\frac{\pi}{2^{2N+2}}.
\]
Thus
\[
\frac{|D(p,r)\cap S_1|}{|D(p,r)|}\geq\frac{|F\cap S_1|}{|D(p,\frac{1}{2^N})|}\geq\frac{\frac{\pi}{2^{2N+4}}-\frac{\pi}{2^{2N+5}}}{\frac{\pi}{2^{2N}}}=\frac{1}{32}.
\]
It is symmetric for $S_2$, $S_3$, $S_4$.  This completes the proof of the lemma. \hfill $\Box$

\end{document}